\documentstyle[12pt]{article}
\input amssym

\newcommand{\bR}{\bf R}

\newcommand{\Mat}{\rm Mat}
\newcommand{\0}{\bf 0}
\newcommand{\1}{\bf 1}
\newcommand{\cle}{\preccurlyeq}
\newcommand{\clt}{\prec}

\begin{document}

\begin{center}

Programming and Computer Software. Vol.~26, No.~5, 2000, pp.275--280\\
Original Russian publication: Programmirovanie. 2000, No.~5, pp.53--62.

\hrule\smallskip
\hrule

\vskip1cm

{\large
Universal Numerical Algorithms and Their Software

Implementation}

\bigskip

{\sc G.L. Litvinov and E.V. Maslova}

International Sophus Lie Center,

 Nagornaya ul. 27--4--72, Moscow 113186 Russia

e-mail: litvinov@islc.msk.su, glitvinov@mail.ru

Received April 10, 2000

\end{center}

\medskip
\begin{abstract}
\noindent
The concept of a universal algorithm is discussed.
Examples of this kind of algorithms are presented.
Software implementations of such algorithms in $C^{++}$-type
languages are discussed together with means that provide
for computations with an arbitrary accuracy. Particular
emphasis is placed on universal algorithms of linear
algebra over semir\-ings.
\end{abstract}

\medskip

\begin{center}
 INTRODUCTION
\end{center}

 Modern achievements in software development
and mathematics make us consider numerical algorithms and their
classification from a new point of view. Conventional numerical
algorithms are oriented to software (or hardware) implementation with
the use of the floating point arithmetic and fixed accuracy. However,
it is often desirable to perform computations with variable (and
arbitrary) accuracy. For this purpose, algorithms are required
that are independent of the accuracy of computations and of a
particular computer representation of numbers. In fact, many
algorithms are not only independent of the computer representation
of numbers, but also of concrete mathematical (algebraic) operations
on data. In this case, operations may be considered as variables.
Such algorithms are implemented by {\it generic programs} based on
the abstract data types technique (abstract data types are
defined by the user, in addition to predefined types of the
language used). The corresponding program tools appeared as
early as in Simula-67, but modern object-oriented languages (like
$C^{++}$, see, e.g. \cite{1,2}) are more convenient for generic programming.

The concept of a generic program was introduced by many authors;
for example, in \cite{3}, such programs were called program schemes.
In this paper, we discuss {\it universal algorithms} implemented as generic
programs and their specific features. This paper is closely related
to papers \cite{4, 5}, in which the concept of a universal algorithm was
defined and software and hardware implementation of such algorithms was
discussed in connection with problems of idempotent mathematics
\cite{4,6}. In this paper, the emphasis is placed on software
 implementation of universal algorithms, computations with arbitrary
accuracy, universal algorithms of linear algebra over semirings, and
their implementation in $C^{++}$.

\begin{center}
 1. UNIVERSAL ALGORITHMS
\end{center}

 Computational algorithms are constructed
on the basis of certain basic operations. Basic operations manipulate
data that describe ``numbers''. These ``numbers'' are elements of a
``numerical domain'', i.e., a mathematical object like the field of
real numbers, the ring of integers, or an idempotent semiring of numbers
(idempotent semirings and their role in idempotent mathematics
are discussed in \cite{4,6} and below in this paper). In every
particular computation, elements of the numerical domains are replaced
by their computer representations, i.e., by elements of certain finite
models of these domains. Examples of models that can be conveniently used
for computer representation of real numbers are provided by various
modifications of floating point arithmetics, approximate arithmetics of
rational numbers \cite{7}, and interval arithmetics. The difference
between mathematical objects (``ideal'' numbers) and their finite
models (computer representations) results in computational (e.g.,
rounding) errors.

An algorithm is called {\it universal} \enskip if it is independent of a
particular numerical domain and (or) of its computer representation.
A typical example of a universal algorithm is computation of the
scalar product $(x,y)$ of two vectors $x=(x_1,\dots,x_n)$ and
$y=(y_1,\dots,y_n)$ by the formula $(x,y)=x_1y_1+\dots+x_ny_n$.
This algorithm (formula) is independent of a particular
domain and its computer implementation, since the formula is
defined for any semiring. It is clear that one algorithm can be
more universal than another. For example, the simplest
rectangular formula provides the most universal algorithm for
numerical integration; indeed, this formula is valid even for
idempotent integration (over any idempotent semiring \cite {4}).
Other quadrature formulas (e.g., combined trapezoid or Simpson
formulas) are independent of computer arithmetics and can be
used (e.g., in the iterative form) for computations with
arbitrary accuracy. In contrast, algorithms based on
Gauss--Jacobi formulas are designed for fixed accuracy computations:
they include constants (coefficients and nodes of these formulas)
defined with fixed accuracy. Certainly, algorithms of this type can
be made more universal by including procedures for computing the
constants; however, this results in an unjustified complication of the
algorithms.

Computer algebra algorithms used in such systems as Mathematica,
Maple, REDUCE, and others are higly universal. Standard
algorithms used in linear algebra can be rewritten in such a way
that they will be valid over any field and complete idempotent
semiring (including semirings of intervals; see \cite{8,9}, where
an interval version of the idempotent linear algebra and the
corresponding universal algorithms are discussed).

As a rule, iterative algorithms (beginning with the successive approximation
method) for solving differential equations (e.g., methods of
Euler, Euler--Cauchy, Runge--Kutta, Adams, a number of important
versions of the difference approximation method, and the like),
methods for calculating elementary and some special functions based on
the expansion in Taylor's series and continuous fractions
(Pad\'e approximations) and others are independent of the computer
representation of numbers.

\clearpage

\begin{center}
2. UNIVERSAL ALGORITHMS AND ACCURACY

OF COMPUTATIONS
\end{center}

Calculations on computers usually are based on a floating-point arithmetic
with a mantissa of a fixed length; i.e., computations are performed
with fixed accuracy. Broadly speaking, with this approach only
the relative rounding error is fixed, which can lead to a drastic
loss of accuracy and invalid results (e.g., when summing series and
subtracting close numbers). On the other hand, this approach provides
rather a high speed of computations. Many important numerical algorithms
are designed to use a floating-point arithmetic (with fixed accuracy)
and ensure the maximum computation speed. However, these algorithms
are not universal. The above mentioned Gauss--Jacobi quadrature formulas,
computation of elementary and special functions on the basis of the
best polynomial or rational approximations or Pad\'e--Chebyshev
approximations, and some others belong to this type. Such algorithms
use nontrivial constants specified with fixed accuracy.

Recently, problems of accuracy, reliability, and authenticity of
computations (including the effect of rounding errors) have come
to the fore; in part, this fact is related to the ever-increasing
performance of computer hardware. When errors in initial data and
rounding errors strongly affect the computation results (ill-posed
problems, analysis of stability of solutions, etc.), it is often useful
to perform computations with improved and variable accuracy. In
particular, the rational arithmetic, in which the rounding error is
specified by the user \cite{7}, can be used for this purpose.
This arithmetic is a useful complement to the interval analysis
\cite{10}. The corresponding computational algorithms must be
universal (in the sense that they must be independent of the computer
representation of numbers).

\begin{center}
4. MATHEMATICS OF SEMIRINGS
\end{center}

A broad class of universal algorithms is related to the concept of
a semiring. We reiterate here the definition of a semiring (see, e.g.,
 \cite{11}). Let $S$ be a set on which associative binary operations
$\oplus$ and $\odot$, called addition and multiplication, respectively,
are defined. We assume that addition is commutative and that multiplication
is distributive over addition; i.e.,
$x\odot (y\oplus z)= (x\odot y) \oplus (x\odot z)$ and
$(x\oplus y)\odot z = (x\odot z) \oplus (y\odot z)$ for all $x, y, z
\in S$. In this case, $S$ is called a {\it semiring}. We assume that
the semiring $S$ contains {\it identity} $\1$ and {\it zero} $\0$; i.e.,
${\1}\odot x = x\odot {\1} = x$ and
${\0}\oplus x = x$, ${\0}\odot x = x\odot{\0} = {\0}$;
in addition, $\0\neq \1$. As is customary, we sometimes omit the
multiplication symbol.

A semiring $S$ is called {\it commutative} if multiplication
$\odot$ is commutative. A semiring $S$ is called {\it idempotent}
if $x\oplus x = x$ for all $x\in S$. If a semiring is a group under
addition, it is called a {\it ring} (in this case, it cannot be
idempotent). If every nonzero element of a commutative ring (semiring)
is invertible under multiplication, this ring (semiring) is called a
{\it field} ({\it semifield}).

The best known and most important examples of semirings are ``numerical''
semirings consisting of real numbers. For example, the set $\bR$
of all real numbers is a field under ordinary arithmetic operations;
i.e., $\oplus = +$, $\odot = \cdot$, ${\0} = 0$, ${\1} = 1$.
The set ${\bR}_{\max} = {\bR} \cup\{-\infty\}$ equipped with operations
$\oplus = \max$ and $\odot = +$ provides an example of an idempotent
semiring (and semifield). Here ${\0} = -\infty$ and ${\1} = 0$.
This semifield is often called the Max-Plus algebra. The semiring
${\bR}_{\min} = {\bR} \cup\{+\infty\}$ equipped with operations
$\oplus = \min$ and $\odot = +$ is isomorphic to the Max-Plus algebra. Here
${\0} = +\infty$ and ${\1} = 0$. Another example is the set
$S^{[a, b]}_{\max, \min}$ consisting of the elements of an interval
$[a,b]$, where $-\infty \leq a < b\leq +\infty$, equipped with
operations $\oplus = \max$ and $\odot = \min$; here ${\0} = a$ and
${\1} = b$. This commutative semiring is not a semifield.

An important example of a  noncommutative semiring is the set
${\Mat}_n(S)$ of all matrices of order $n\times n$ with  elements
from a commutative semiring $S$ with ordinary standard operations.
The {\it sum} of matrices $A=(a_{ij})$ and $B=(b_{ij})$ is the
matrix $A\oplus B = (a_{ij}\oplus b_{ij})$, and the {\it product}
of these matrices is the matrix
$AB =(\bigoplus^n_{k=1} a_{ik}\odot b_{kj})$, where $i, j = 1, \dots, n$.
Operations on rectangular matrices can be defined similarly. Zero
$O$ and identity $I$ in ${\Mat}_n (S)$ are defined in the conventional way.
If the semiring $S$ is idempotent, then ${\Mat}_n (S)$ is also idempotent.
Many other important examples can be found in \cite{3} -- \cite{6},
\cite{8}, \cite{9}, \cite{11}.

On any idempotent semiring, a  {\it canonical partial order} $\cle$
is defined by the following rule: $x\cle y$ is equivalent to
$x\oplus y = y$. Moreover, $x\oplus y = \sup\{ x, y\}$ with respect
to the canonical order. The canonical order is
compatible with the semiring addition and
multiplication in the common way. For the semirings ${\bR}_{\max}$
and $S^{[a,b]}_{\max, \min}$, the canonical order coincides with the
standard order $\leq$ defined on the set of real numbers; for the
semiring $\bR_{\min}$, it is inverse to the standard order.

There exists a (heuristic) correspondence between important,
useful, and interesting constructs and results of traditional
mathematics over fields and similar constructs and results of
idempotent mathematics (i.e., mathematics over idempotent semirings).
This idempotent correspondence principle is closely related to the
Bohr correspondence principle in quantum mechanics. Traditional
mathematics can be considered as a ``quantum'' theory and idempotent
mathematics as its  ``classical'' analogue (see \cite{4}).
Consistent application of the idempotent correspondence principle leads to
various and surprising results, including a methodology for
constructing universal algorithms and patenting computer devices
 \cite{4}, \cite{5}.

The fundamental equations in quantum theory are linear
(superposition principle). There is an idempotent version of the
superposition principle \cite{6}): the Hamilton--Jacobi equation,
i.e., the basic (nonlinear) equation of classical mechanics, can
be considered as linear over the semiring ${\bR}_{\min}$; various
modifications of the Bellman equation, i.e., the basic equation of
optimization theory, are also linear over appropriate idempotent
semirings. For example, the finite-dimensional time-independent
Bellman equation can be written as
\begin{equation}
X=A\odot X\oplus B,
\end{equation}
where $A$ is a square matrix with elements from an idempotent semiring
$S$ and $X$ and $B$ are vectors (or matrices) with elements from $S$.
The solution $X$ is found from (1) when $A$ and $B$ are given.

In particular, standard problems in dynamic programming correspond
to the case $\bR_{\max}$, and the well-known shortest path problem corresponds
to $S={\bR}_{\max}$. It is shown in \cite{12} that the principal
optimization algorithms for finite graphs correspond to standard
methods for solving systems of linear equations of form (1) over semirings.
The Bellman algorithm for the shortest path problem corresponds to
a semiring version of the Jacobi method; the Ford algorithm
corresponds to the Gauss--Seidel iterative method; and so on. These
algorithms are universal and may be used for solving linear algebra
problems over a broad class of semirings that includes all idempotent
semirings and all fields.

Idempotent analogues of standard numerical algorithms are very
important and can be used systematically for solving, for example,
optimization problems. Linear algebra algorithms are of prime
importance, since standard infinite-dimensional linear problems
over semirings can be reduced to finite-dimensional (or finite)
approximations, and nonlinear algorithms can often be approximated
by linear ones.

We note that the available methods used for parallelizing linear algebra
algorithms can be applied to their semiring analogues.

\begin{center}
5. UNIVERSAL LINEAR ALGEBRA

ALGORITHMS OVER SEMIRINGS
\end{center}

The most important linear algebra problem is solving the system
of linear equations
\begin{equation}
AX = B,
\end{equation}
where $A$ is a matrix with elements from the basic field and $X$ and
$B$ are vectors (or matrices) with elements from the same field.
It is required to find $X$ if $A$ and $B$ are given. If $A$ in (2)
is not the identity matrix $I$, then
system (2) can be written in form (1), i.e.,
$$
X = AX + B.\eqno{(1')}
$$
It is well known that form (1) or ($1'$) is convenient for using the
successive approximation method. Applying this method with the initial
approximation $X_0=0$, we obtain the solution
\begin{equation}
X = A^*B,
\end{equation}
where
\begin{equation}
A^* = I+A+A^2+\dots + A^n+\dots
\end{equation}
On the other hand, it is clear that
\begin{equation}
A^* = (I-A)^{-1},
\end{equation}
if the matrix $I-A$ is invertible. The inverse matrix $(I-A)^{-1}$
can be considered as a regularized sum of the formal series (4).

The above considerations can be extended to a broad class of
semirings. The unary operation $A\mapsto A^*$ in ${\Mat}_n(S)$
is defined (partially) if a unary (partial) operation $x\mapsto x^*$,
called {\it closure}, is defined on the semiring $S$ such that
the identity
\begin{equation}
x^* = {\1}\oplus (x^*\odot x) = {\1}\oplus (x\odot x^*)
\end{equation}
holds true if $x^*$ is defined. It follows from (6) that
$$
x^* = {\1}\oplus x\oplus x^2\oplus\dots\oplus x^n\oplus x^*x^{n+1}
$$
for any positive integer $n$; thus, $x^*$ can be considered as
a regularized sum of the formal series
$$
x^* = {\1}\oplus x\oplus x^2\oplus\dots\oplus x^n\oplus\dots
$$

If $S$ is a field, then, by definition, $x^*=({\1}-x)^{-1}$ for any
$x\neq \1$. If $S$ is an idempotent semiring, then, by definition
\begin{equation}
x^* = {\1}\oplus x\oplus x^2\oplus\dots= \sup\{{\1}, x, x^2,\dots,x^n
\dots\},
\end{equation}
if this supremum (with respect to the canonical order $\cle$)
exists. In this case, $x^*={\1}$ if $x\cle {\1}$. Therefore, $x^*={\1}$
in the semiring $S^{[a,b]}$ for all $x$. For the semifield ${\bR}_{\max}$
the closure operator $x\mapsto x^*$ is not defined for ${\1}\clt x$
(however, ${\bR}_{\max}$ can be supplemented by $+\infty$, which
turns this semifield into a semiring; in this case, $x^*=+\infty$ for
${\1}\clt x$). It is clear that, for $x\cle {\1}$, $x^*={\1}$ in
${\bR}_{\max}$, as well as in other idempotent semirings. These
examples show that the closure $x^*$ of $x$ is often calculated very
simply for idempotent semirings.

The closure operation for matrix semirings ${\Mat}_n(S)$ can be defined
and computed in terms of the closure operation for $S$; some
methods are described in \cite{3, 5, 6, 11, 12}. One such method is
described below ($LDM$-factorization). The closure operation
$A\mapsto A^*$ in ${\Mat}_n(S)$ satisfies identity (6), which implies
that if $A^*$ is defined, then $X=A^*B= A^*\odot B$ is the solution
to the matrix equation (1).

Consider a nontrivial universal algorithm applicable to matrices over
semirings with the closure operation defined.

\begin{center}
{\it Example: Semiring $LDM$-Factorization}
\end{center}

Factorization of a matrix into the product $A = LDM$, where $L$ and $M$
are lower and upper triangular matrices with a unit diagonal,
respectively, and $D$ is a diagonal matrix, is used for solving
matrix equations $AX = B$ \cite{13}. We construct a similar
decomposition for the Bellman equation $X = AX \oplus B$.

For the case $AX = B$, the decomposition $A = LDM$ induces the following
decomposition of the initial equation:
\begin{equation}
   LZ = B, \qquad DY = Z, \qquad MX = Y.
\end{equation}
Hence, we have
\begin{equation}
   A^{-1} = M^{-1}D^{-1}L^{-1},
\label{AULinv}
\end{equation}
if $A$ is invertible. In essence, it is sufficient to find the
matrices $L$, $D$ and $M$, since the linear system (8) is easily
solved by a combination of the forward substitution for $Z$, the
trivial inversion of a diagonal matrix for $Y$, and the back
substitution for $X$.

Using (8) as a pattern, we can write
\begin{equation}
   Z = LZ \oplus B, \qquad Y = DY \oplus Z, \qquad X = MX \oplus Y.
\label{LDM}
\end{equation}
Then
\begin{equation}
   A^* = M^*D^*L^*.
\label{AMDLstar}
\end{equation}
A triple $(L,D,M)$ consisting of a lower triangular, diagonal, and
upper triangular matrices is called an $LDM$-{\it factorization} of a
matrix $A$ if relations (10) and (11) are satisfied. We note that
in this case, the principal diagonals of $L$ and $M$ are zero.

The modification of the notion of $LDM$-factorization used in matrix
analysis for the equation $AX=B$ is constructed by analogy with the
construct suggested by Carr\'e in \cite{12} for $LU$-factorization.

We stress that the algorithm described below can be applied to matrix
computations over any semiring under the condition that the unary
operation $a\mapsto a^*$ is applicable every time it is encountered
in the computational process. Indeed, when constructing the
algorithm, we use only the basic semiring operations of addition
$\oplus$ and multiplication $\odot$ and the properties of
associativity, commutativity of addition, and distributivity of
multiplication over addition.

If $A$ is a symmetric matrix over a semiring with a commutative
multiplication, the amount of computations can be halved, since
$M$ and $L$ go into each other under transposition.

We begin with the case of a triangular matrix $A = L$ (or $A = M$).
Then, finding $X$ is reduced to the forward (or back) substitution.

\begin{center}
{\it Forward substitution}
\end{center}

 We are given:
\begin{itemize}
\item $L = \|l^i_j\|^n_{i,j = 1}$, where $l^i_j = \0$ for $i \le j$
(a lower triangular matrix with a zero diagonal);
\item $B = \|b^i\|^n_{i = 1}$.
\end{itemize}

It is required to find the solution $X = \|x^i\|^n_{i = 1}$ to the
equation $X = LX \oplus B$. The program fragment solving this problem is as
follows.

\begin{tabbing}
   \qquad\=\qquad\=\kill
   for $i = 1$ to $n$ do\\*
   \{\> $x^i := b^i$;\\
   \> for $j = 1$ to $i - 1$ do\\*
   \>\> $x^i := x^i \oplus (l^i_j \odot x^j)$;\, \}\\
\end{tabbing}

\begin{center}
{\it Back substitution}
\end{center}

We are given
\begin{itemize}
\item $M = \|m^i_j\|^n_{i,j = 1}$, where $m^i_j = \0$ for $i \ge j$ (an
upper triangular matrix with a zero diagonal);
\item $B = \|b^i\|^n_{i = 1}$.
\end{itemize}

It is required to find the solution $X = \|x^i\|^n_{i = 1}$ to the
equation $X = MX \oplus B$. The program fragment solving this problem
is as follows.

\begin{tabbing}
   \qquad\=\qquad\=\kill
   for $i = n$ to 1 step $-1$ do\\*
   \{\> $x^i :=  b^i$;\\
   \> for $j = n$ to $i + 1$ step $-1$ do\\*
   \>\> $x^i :=  x^i \oplus (m^i_j \odot x^i)$;\, \}\\
\end{tabbing}

Both algorithms require $(n^2 - n) / 2$ operations $\oplus$ and $\odot$.

\begin{center}
{\it Closure of a diagonal matrix}
\end{center}

We are given
\begin{itemize}
\item $D = {\rm{diag}}(d_1, \ldots, d_n)$;
\item $B = \|b^i\|^n_{i = 1}$.
\end{itemize}

It is required to find the solution $X = \|x^i\|^n_{i = 1}$ to the
equation $X = DX \oplus B$. The program fragment solving this problem
is as follows.

\begin{tabbing}
   \qquad\=\qquad\=\kill
   for $i = 1$ to $n$ do\\*
   \> $x^i :=  (d_i)^* \odot b^i$;\\
\end{tabbing}

This algorithm requires $n$ operations $*$ and $n$ multiplications $\odot$.

\begin{center}
{\it General case}
\end{center}

We are given

\begin{itemize}
\item $L = \|l^i_j\|^n_{i,j = 1}$, where $l^i_j = \0$ if $i \le j$;
\item $D = {\rm{diag}}(d_1, \ldots, d_n)$;
\item $M = \|m^i_j\|^n_{i,j = 1}$, where $m^i_j = \0$ if $i \ge j$;
\item $B = \|b^i\|^n_{i = 1}$.
\end{itemize}

It is required to find the solution $X = \|x^i\|^n_{i = 1}$ to the
equation $X = AX \oplus B$, where $L$, $D$, and $M$ form the
$LDM$-factorization of $A$. The program fragment solving this problem
is as follows.

\begin{tabbing}
        {\sc {FORWARD SUBSTITUTION}}\\*
   for $i = 1$ to $n$ do\\*
   \{\, $x^i :=  b^i$;\\*
   \, for $j = 1$ to $i - 1$ do\\*
   \,\, $x^i :=  x^i \oplus (l^i_j \odot x^j)$;\, \}\\
	\sc{CLOSURE OF A DIAGONAL MATRIX}\\*
   for $i = 1$ to $n$ do\\*
   \, $x^i :=  (d_i)^* \odot b^i$;\\
	\sc{BACK SUBSTITUTION}\\*
   for $i = n$ to 1 step $-1$ do\\*
   \{\, for $j = n$ to $i + 1$ step $-1$ do\\*
   \,\, $x^i :=  x^i \oplus (m^i_j \odot x^j)$;\, \}\\
\end{tabbing}

Note that $x^i$ is not initialized in the course of the back substitution.
The algorithm requires $n^2 - n$ operations $\oplus$, $n^2$ operations
$\odot$, and $n$ operations~$*$.

\begin{center}
{\it LDM-factorization}
\end{center}

We are given
\begin{itemize}
\item $A = \|a^i_j\|^n_{i,j = 1}$.
\end{itemize}

It is required to find the $LDM$-factorization of $A$:
$L = \|l^i_j\|^n_{i,j = 1}$, $D ={\rm{diag}}(d_1, \ldots, d_n)$, and
$M = \|m^i_j\|^n_{i,j = 1}$, where $l^i_j = \0$ if $i \le j$, and
$m^i_j = \0$ if $i \ge j$.

The program uses the following internal variables:
\begin{itemize}
\item $C = \|c^i_j\|^n_{i,j = 1}$;
\item $V = \|v^i\|^n_{i = 1}$;
\item $d$.
\end{itemize}

\begin{tabbing}
   \qquad\=\qquad\=\qquad\=\kill
   \sc{INITIALISATION}\\*
	for $i = 1$ to $n$ do\\*
	\> for $j = 1$ to $n$ do\\*
	\>\> $c^i_j = a^i_j$;\\
	\sc{MAIN LOOP}\\*
	for $j = 1$ to $n$ do\\*
	\{\> for $i = 1$ to $j$ do\\*
	\>\> $v^i :=  a^i_j$;\\
	\> for $k = 1$ to $j - 1$ do\\*
	\>\> for $i = k + 1$ to $j$ do\\*
	\>\>\> $v^i :=  v^i \oplus (a^i_k \odot v^k)$;\\
	\> for $i = 1$ to $j - 1$ do\\*
	\>\> $a^i_j :=  (a^i_i)^* \odot v^i$;\\
	\> $a^j_j :=  v^j$;\\
	\> for $k = 1$ to $j - 1$ do\\*
	\>\> for $i = j + 1$ to $n$ do\\*
	\>\>\> $a^i_j :=  a^i_j \oplus (a^i_k \odot v^k)$;\\
	\> $d = (v^j)^*$;\\
	\> for $i = j + 1$ to $n$ do\\*
	\>\> $a^i_j :=  a^i_j \odot d$;\, \}\\
\end{tabbing}

This algorithm requires $(2n^3 - 3n^2 + n) /6$ operations $\oplus$, $(2n^3 +
3n^2 -5n) / 6$ operations $\odot$, and $n(n + 1) / 2$ operations $*$.
After its completion, the matrices $L$, $D$, and $M$ are contained,
respectively, in the lower triangle, on the diagonal, and in the upper
triangle of the matrix $C$. In the case when $A$ is symmetric about the
principal diagonal and the semiring over which the matrix is defined
is commutative, the algorithm can be modified in such a way that the
number of operations is reduced approximately by a factor of two. For
details see \cite{13}.

Other examples can be found in \cite{3}, \cite{11} -- \cite{15}.

\begin{center}
6. SOFTWARE IMPLEMENTATION

OF UNIVERSAL ALGORITHMS
\end{center}

Object-oriented languages (e.g., $C^{++}$ and Java) and programming
systems that allow abstract data types to be defined provide
convenient means for the software implementation of universal
algorithms. In this case, program units can operate with abstract
(and variable) operations and data types. Specific values of
operations are determined by the input data types, these operations
(and data types) are implemented by additional program units.
Recently, this type of programming technique has been dubbed
generic programming (see, e.g., \cite{1,2}). To help automate the
generic programming, the so-called Standard Template Library (STL)
was developed in the framework of $C^{++}$ \cite{2,16}. However,
high-level tools, such as STL, possess both obvious advantages
and some disadvantages and must be used with caution.

Using the generic programming technique, a program package was
developed in $C^{++}$ for solving problems in linear algebra over
fields and semirings (for various computer implementation of the
corresponding numeric domains) and optimization problems on
graphs. A hierarchy of abstract data types for basic numeric fields,
rings, semifields, and semirings was developed for various computer
representations. In particular, various versions of the rational
arithmetic \cite{7} can be used and computations can be performed with
any given accuracy. Solving systems of linear Bellman equations over
idempotent semirings (by various methods), standard optimization
problems on graphs can be solved (the dynamic programming problem,
shortest path problem, widest path problem, etc.), including
interval versions of those problems \cite{8}, \cite{9}. The system
provides a basis for a more powerful program package based on universal
algorithms \cite{5}. This system will be described in detail in
subsequent publications.

\begin{center}
ACKNOWLEDGMENTS
\end{center}

This work was supported by the Russian Foundation for Basic Research,
project no. 99--01--01 198. We are grateful to A. Ya. Rodionov and A. N.
Sobolevskii for useful advice and help.

\end{document}